\documentstyle[12pt]{article}
\newcommand{\const}{\mathop{\rm const}\limits}

\newcommand{\Div}{\mathop{\rm Div}\limits}

\begin{document}

\begin{center}

{\bf  QUANTITATIVE \\
 LOWER BOUND FOR LIFESPAN  FOR SOLUTION \\
OF NAIVER-STOKES EQUATIONS.} \par

\vspace{4mm}

 $ {\bf E.Ostrovsky^a, \ \ L.Sirota^b } $ \\

\vspace{4mm}

$ ^a $ Corresponding Author. Department of Mathematics and computer science, Bar-Ilan University, 84105, Ramat Gan, Israel.\\
\end{center}
E - mail: \ galo@list.ru \  eugostrovsky@list.ru\\
\begin{center}
$ ^b $  Department of Mathematics and computer science. Bar-Ilan University,
84105, Ramat Gan, Israel.\\

E - mail: \ sirota3@bezeqint.net\\

\vspace{5mm}
                    {\bf Abstract.}\\

 \end{center}

 \vspace{4mm}

 We find a simple quantitative  {\it lower} bound for lifespan of solution  of the multidimensional initial value
problem for the Navier-Stokes equations in whole space when the initial function belongs to the correspondent Lebesgue - Riesz space,
and give some a priory estimations for solution in some rearrangement invariant spaces. \par

\vspace{4mm}

 \par

\vspace{4mm}

  \par

\vspace{4mm}

\section{Notations. Statement of problem.}

\vspace{3mm}

{\bf Statement of problem.} \par

\vspace{3mm}

We  consider in this article the  initial value problem  for the multivariate Navier-Stokes (NS) equations

$$
\partial{u}_t - \Delta u + (u \cdot \nabla)u = \nabla P, \ x \in R^d, \ d \ge 3,   \ t > 0; \eqno(1.1)
$$

$$
\Div (u) = 0, \ x \in R^d, \ t > 0; \eqno(1.2)
$$

$$
u(x,0) = a(x), \ x \in R^d. \eqno(1.3)
$$
 Here as ordinary

$$
x = (x_1,x_2,\ldots,x_k,\ldots,x_d) \in R^d ; \  ||x||: = \sqrt{\sum_{j=1}^d x_j^2}.
$$
and
$$
u = u(t) = u(t,\cdot)  = u(x,t)  = \{ u_1(x,t), u_2(x,t), \ldots,  u_d(x,t) \}
$$
denotes the (vector) velocity of fluid in the point $ x $ at the time $  t,  \ P $   is
represents the pressure. \par
 Equally:

 $$
\partial{u_i}/\partial t  = \sum_{j=1}^d \partial^2_{x_j} u_i - \sum_{j=1}^d u_j \partial_{x_j} u_i +
+ \partial_{x_i} P,
 $$

$$
\sum_{j=1}^d \partial_{x_j} u_j = 0, \ u(x,0) = a(x),
$$

$$
\Div u = \Div \vec{u} = \Div \{ u_1, u_2, \ldots, u_d \} = \sum_{k=1}^d \frac{\partial u_k }{\partial x_k} = 0
$$
in the sense of distributional derivatives.\par
 As long as

 $$
 P = \sum \sum_{j,k = 1}^d R_j R_k (u_j \cdot u_k),
 $$
where $ R_k = R_k^{(d)} $ is the $ k^{th} \ d \- $ dimensional Riesz transform:

$$
 \ R_k^{(d)}[f](x) = c(d) \lim_{\epsilon \to 0+} \int_{ ||y|| > \epsilon} ||y||^{-d} \Omega_k(y) \ f(x-y) \ dy,
$$

$$
c(d) = -\frac{\pi^{(d+1)/2}}{\Gamma \left( \frac{d+1}{2} \right) }, \ \Omega_k(x)= x_k /||x||,
$$
the system  (1.1) - (1.3) may be rewritten as follows:

$$
\partial{u}_t = \Delta u + (u \cdot \nabla)u  +   Q \cdot \nabla \cdot (u \otimes u), \ x \in R^d, \ t > 0; \eqno(1.4)
$$

$$
\Div (u) = 0, \ x \in R^d, \ t > 0; \eqno(1.5)
$$

$$
u(x,0) = a(x), \ x \in R^d, \eqno(1.6)
$$
where   $ Q $ is multidimensional  Helmholtz-Weyl projection
operator, i.e., the $ d \times d  $ matrix pseudo-differential operator in $ R^d $ with the matrix symbol

$$
a_{i,j}(\xi) = \delta_{i,j} - \xi_i \xi_j /||\xi||^2, \hspace{5mm} \delta_{i,j} = 1, i = j; \delta_{i,j} = 0, \ i \ne j.
$$

\vspace{4mm}
 {\it We will understand henceforth as a capacity of the solution (1.4) - (1.6) the vector - function  } $ u = \vec{u} =
 \{ u_1(x,t),  u_2(x,t), \ldots,  u_d(x,t) \} $ {\it  the so-called mild solution,} see \cite{Muira1}.  \par
\vspace{4mm}

 Namely, the vector- function  $ u = u(t)  $ satisfies almost everywhere in the time $ t $ the following {\it non-linear
 integral - differential equation: }

$$
u(t) = e^{t \Delta} a + \int_0^t e^{(t-s)\Delta } [ (u \cdot \nabla)u(s)  +   Q \cdot \nabla \cdot (u \otimes u)(s) ] ds \stackrel{def}{=}
$$

$$
e^{t \Delta} a  + G [u](t) \stackrel{def}{=} u_0(x,t) + G [u](t), \eqno(1.7)
$$
 the operator $  \exp(t \Delta) $ is the classical convolution  integral operator with heat kernel: \par

$$
u_0(x,t) :=  e^{t \Delta}[ a](x,t) = w_t(x)*a(x),
$$
where

$$
G(u)  \stackrel{def}{=}  F(u,u) = F(u),  \hspace{4mm}  F(u,v):= \int_0^t[ (u \cdot \nabla)v(s)  +   Q \cdot \nabla \cdot (u \otimes v)(s) ] \ ds,
$$

$$
w_t(x) \stackrel{def}{=} (4 \pi t)^{-d/2}  \exp \left( - \frac{|| x ||^2}{4 t} \right) \eqno(1.8)
$$
The convolution between two functions   $ r = f(t), \ g(t)  $ defined on the set $ R_+ = (0,\infty)  $ is defined as usually

$$
f \odot g(t) = \int_0^t f(t-s) \ g(s) \ ds  = g \odot f(t)
$$
("time=wise" convolution)
and between two, of course, measurable  vector - functions  $ u(x), v(x) $ defined on the whole space $ x \in R^d $

$$
u*v(x) = \int_{R^d} u(x-y) \ v(y) \ dy,
$$
("space wise" and  coordinate-wise convolution).
  The authors hope that this notations does not follow the confusion. \par

\vspace{3mm}

   More results about the existence, uniqueness, numerical methods, and a priory estimates in the different Banach function spaces:
 Lebesgue-Riesz $  L_p, $ Morrey, Besov for this solutions see, e.g. in  \cite{Cui1}- \cite{Zhang1}. The first and besides famous
 result belong to J.Leray  \cite{Leray1}; it is established there in particular the {\it global in time}
 solvability and  uniqueness of NS system in the space $ L_2(R^d)  $ and was obtained  a very interest a priory estimate for solution.\par
  The immediate predecessor for offered article is  the article of T.Kato \cite{Kato1}; in this article was considered the case
  $  a(\cdot) \in L_d(R^d).$  See also celebrate works of H. Fujita  and T.Kato \cite{Fujita1},
   Y.Giga  \cite{Giga1} - \cite{Giga4}, T.Kato  \cite{Kato2} etc.\par

  T.Kato in \cite{Kato1} proved in particular that if the initial function $ a = a(x)  $ belongs to the space $  L_d(R^d) $
(in our notations), then there exists a positive time value $ T > 0  $ (lifespan of solution) such that the solution of NS system
$ u = u(x,t) $ there exists for $ t \in (0,T), $  is smooth and satisfy some a priory integral estimates.\par
 Furthermore,  if the norm $ ||a||L_d(R^d)  $ is sufficiently small, then $  T = \infty, $ i.e. the solution $ u = u(x,t) $
is global.\par
 The {\it upper} estimate for the value $ T, $ conditions for finite-time blow-up
 and asymptotical behavior of solution as $ t \to T - 0 $ see in the articles \cite{Ball1}, \cite{Benameur1},
\cite{Chae1},  \cite{Cui1}, \cite{Gallagher1}, \cite{Montgomery1},  \cite{Seregin1}, \cite{Seregin2} etc.\par

\vspace{4mm}

{\bf  Our purpose in this report is to obtain the quantitative simple lower estimates for the lifespan of solution } $  T. $ \par

\vspace{4mm}

 For the chemotaxis equations under some  additional conditions this estimate is obtained in
\cite{Hatami1}. \par

\vspace{4mm}

 In detail, we understand  as a solution $ u = u(x,t), \ t \in (0,T)  $  together with T.Kato \cite{Kato1}
  the mild solution of NS equations such that for all the values $  q \ge d  $

$$
t^{ (1-d/q)/2 } \cdot u  \in BC([0,T], L_q^0),   \hspace{5mm} t^{ 1-d/2q } \cdot \nabla u  \in BC([0,T], L_q^0). \eqno(1.9)
$$

The critical value of the variable $  T,  $  more exactly, for its lower estimate does not depend on the variable $  q. $ \par

 Here $  T \in (0,\infty); $ the case $ T = \infty $ implies the global (in time) solution. \par

 The space  $ BC([0,T], L_q^0) $ consists by definition  on all the functions $ v = v(x,t), \ x \in R^d, \ t \in [0, T] $
with zero divergence and finite norm

$$
|| v(\cdot, \cdot) ||BC([0,T], L_q^0) := \sup_{t \in [0,T] } ||v(\cdot,t)||_q.  \eqno(1.10)
$$

 Will  be also presumed for all the functions from the space $ BC([0,T], L_q^0) $  the  { \it continuity} on the time $  t $
in the $ L_q, \ q \ge d $  sense:

 $$
 \lim_{s \to t} || v(\cdot,t) - v(\cdot,s)||_q = 0, \ t \in [0,T].
 $$

 At the value $ t = 0+ $

$$
 \lim_{t \to 0+} ||u(\cdot,t) - u(\cdot,0)||_q = 0, \ q \ge d.
$$

 For instance, the initial condition $  a(x) = u(x,0) $ will be understood as follows:

 $$
 \lim_{t \to 0+} ||u(\cdot,t) - a(\cdot)||_d = 0. \eqno(1.11)
 $$

 \vspace{3mm}

 This estimates allow us to  establish some new properties of solution and develop numerical methods. \par

 Note that this statement of problem appeared in \cite{Ostrovsky111}. \par

\vspace{3mm}

\section{ Some Notations, with Clarification.} \par

\vspace{3mm}

 As ordinary, for the measurable function $ x \to u(x), \ x \in R^d $

 $$
||u||_p = \left[ \int_{R^d} |u(x)|^p \ dx  \right]^{1/p}. \eqno(2.1)
 $$

{\it  Multidimensional case.} \par
 Let $ u = \vec {u} = \{ u_1(x), u_2(x), \ldots, u_d(x)  \}  $ be measurable vector - function: $ u_k: R^d \to R. $ We can define
as ordinary the $ L_p, \ p \ge 1 $  norm of the function $ u $ by the following way:

$$
||u||_p := \max_{k=1,2, \ldots,d} || u_k||_p, \ p \ge 1.
$$

\vspace{5mm}

 Define also

$$
K_S(d,p) := \pi^{-1/2} \ d^{-1/p} \ \left( \frac{p-1}{d-p} \right)^{ (p-1)/p } \
\left\{ \frac{\Gamma(1+d/2) \ \Gamma(d)}{\Gamma(d/p) \ \Gamma(1+d - d/p) } \right\}^{1/d}. \eqno(2.4)
$$

 The function $ K_S(d,p) $ is the optimal (i.e. minimal) value in the famous Sobolev's inequality

 $$
 ||\phi||_r \le K_S(d,q) \ || \nabla \phi  ||_q, \ 1 \le q < d, \ \frac{1}{r} = \frac{1}{q} -  \frac{1}{d}, \ r \ge 1, \eqno(2.5)
 $$
see  Bliss \cite{Bliss1}, (1930); Talenti, \cite{Talenti1}, (1995).\par

\vspace{5mm}

$$
\tilde{\omega}(d) := \frac{4 \pi^{d/2-1}}{\Gamma(d/2)},\ \omega(d) := \frac{2 \pi^{d/2}}{\Gamma(d/2)},
$$

\vspace{3mm}

$$
c(d) = -\frac{\pi^{(d+1)/2}}{\Gamma \left( \frac{d+1}{2} \right) }. \hspace{6mm} \Omega_k(x)= x_k  / ||x||.
$$

\vspace{3mm}

$$
x = (x_1,x_2,\ldots,x_k,\ldots,x_d) \in R^d  \ \Rightarrow ||x|| = \sqrt{\sum_{j=1}^d x_j^2}.
$$

\vspace{3mm}

 The explicit view for Riesz's transform has a view
$$
 R_k[f](x)= R_k^{(d)}[f](x) = c(d) \lim_{\epsilon \to 0+} \int_{ ||y|| > \epsilon} ||y||^{-d} \Omega_k(y) \ f(x-y) \ dy.\eqno(2.6)
$$

 The ultimate result in this direction belongs to T.Iwaniec  and  G.Martin  \cite{Iwaniec3}:  the  {\it upper estimate} value
$ ||R_k||(L_p \to L_p) $ does not dependent on the dimension $ d $ and coincides with the Pichorides constant:

$$
K_R(p):=  ||R_k||(L_p \to L_p) = \cot \left( \frac{\pi}{2p^*}  \right), \ p^* = \max(p, \ p/(p-1) ), \ p > 1.\eqno(2.7)
$$
 For instance, $ K_R(3) = \sqrt{3}. $ \par
 T.Iwaniec  and  G.Martin considered  also the vectorial Riesz transform. \par
 See for additional information  \cite{Banuelos1},  \cite{Stein1},  chapter 2, section 4; \cite{Taylor1}, chapter 3. \par

\vspace{4mm}

 We will use the famous Young inequality for the (measurable) functions  $ f,g: R^d \to R:$

 $$
 || f*g||_r  \le K_{B L}(p,q) \ ||f||_p \ ||g||_q, \ 1/r +1 = 1/p + 1/q, \ p,q,r > 1,
 $$
where $ d $ is a dimension of arguments of a functions $   f,g $  and
the optimal value of "constant"   $ K_{B L}(p,q) $ was obtained by
 H.J.Brascamp and E.H. Lieb \cite{Brascamp1}:

 $$
  K_{B L}(d;p,q) = K_{B L}(p,q) = \left[p^{1/p} \ s^{-1/s} \ q^{1/q} \ t^{-1/t} \ r^{1/r} \ z^{-1/z} \right]^{d/2}, \eqno(2.8)
 $$

$$
s = p/(p-1), \ t = q/(q-1), \ z = r/(r-1).
$$
 Note that  $ K_{B L}(d; p,q) \le 1. $ \par

\vspace{3mm}

  Let us denote

  $$
  M(d,r) = ||w_1(\cdot)||_r, \ r \ge 1.
  $$

 We deduce by direct computation

 $$
  M(d,r) =  (4 \pi)^{- d/2}  \left[ \int_{R^d} e^{ - r ||x||^2/(4t)  } dx \right]^{1/r} =  2^{d/r} \ \pi^{-d(1-1/r)/2 } \ r^{-d/2r}.
 $$

 Therefore
 $$
 ||w_t(\cdot)||_r = t^{-d( 1 - 1/r)/2 } M(d,r) = t^{-d( 1 - 1/r)/2 } \ 2^{d/r} \ \pi^{-d(1-1/r)/2 } \ r^{-d/2r}. \eqno(2.9)
 $$

 Note that the value $  M(d,r)  $ allows a very simple estimates: $ M(d,r)  < 2^d. $

\vspace{4mm}
 Let us denote  $ K_0 =  K_0(d,\delta) = K_9(a; d.\delta) := $

 \vspace{4mm}

 $$
 \sup_{t > 0} \left[ t^{(1-\delta)/2 } \ ||u_0(t)||_{d/\delta} \right] =
 \sup_{t > 0} \left[t^{(1-\delta)/2 } \  || e^{t \Delta}  a||_{d/\delta} \right]  =
 $$

 $$
 \sup_{t > 0} \left[ t^{(1-\delta)/2} \ || w_t *a  ||_{d/\delta} \right].
 $$

 We conclude as a consequence using Young inequality:
  if $ t > 0, \  a \in L_d(R^d), \ \delta = \const \in (0,1), $ then

$$
K_0(d,\delta)  \le K_{BL}(d; d, d/(d-(1+ \delta))) \cdot M(d, d/(d-(1 + \delta )) )  \cdot ||a||_d.
\eqno(2.10)
$$

 \vspace{4mm}

 We get analogously denoting  $ K'_0= K'_0(d) = K'_0(a;d) :=    $

 $$
 \sup_{t > 0} \left[ t^{1/2} \cdot || \nabla u_0(t) ||_d \right] =  \sup_{t > 0} \left[ t^{1/2} \ ||  \nabla e^{t \Delta} a||_d \right]   =
 \sup_{t > 0} \left[ t^{1/2} \ || (\nabla w_t) * a ||_d \right]:
 $$

 $$
K'_0 \le 0.5 \cdot K_{BL}(d; d,d) \cdot M \left(d, \frac{d^2}{d - 1} \right)  \cdot ||a||_d. \eqno(2.11)
 $$

\vspace{5mm}

 Notice \cite{Kato1} that for the operator (non-linear) $  G $ are true the following estimates:

 $$
 || Gu ||_{ d/\gamma} \le K_R(d/\alpha) \ K_R(d/\beta)  \
 $$

 $$
  \int_0^t (t-s)^{- (\alpha + \beta - \gamma)/2 } \ ||u(s)||_{d/\alpha} \ || \nabla u(s) ||_{d/\beta} \ ds = \eqno(2.12)
 $$

$$
K_R(d/\alpha) \ K_R(d/\beta)  \
\left\{ t^{- (\alpha + \beta - \gamma)/2 } \odot \left[  \ ||u(t)||_{d/\alpha} \ || \nabla u(t) ||_{d/\beta} \right] \right\};
$$

\vspace{5mm}

 $$
 || \nabla Gu ||_{ d/\gamma} \le K_R(d/\alpha) \ K_R(d/\beta) \times
 $$

 $$
  \int_0^t (t-s) ^{- (1+\alpha + \beta - \gamma)/2}  \ ||u(s)||_{d/\alpha} \ || \nabla u(s) ||_{d/\beta} \ ds=  \eqno(2.13)
 $$

$$
K_R(d/\alpha) \ K_R(d/\beta) \
\left\{ t^{- (1+\alpha + \beta - \gamma)/2} \odot \left[ ||u(t)||_{d/\alpha} \ || \nabla u(t) ||_{d/\beta} \right] \right\},
$$

$$
\alpha, \beta, \gamma > 0,  \gamma \le \alpha + \beta < d.
$$

\vspace{5mm}

 Another useful inequalities: $ || e^{t \Delta} v||_q \le $

 $$
  K_{BL}  \left(d; \frac{pq}{pq + p - q},p  \right) \cdot M \left(d, \frac{pq}{pq + p - q}    \right)   \cdot ||v||_p
  \cdot t^{-(d/p - d/q )/2   }, \ 1 < p \le q < \infty; \eqno(2.14)
 $$

$$
 || \nabla  e^{t \Delta} v||_q \le 0.5 \cdot  K_{BL}\left(d; \frac{pq}{pq + p - q},p  \right) \cdot
$$

$$
  M \left(d; d+\frac{pq}{pq + p - q}  \right) \cdot ||v||_p  \cdot t^{-(1+ d/p - d/q )/2   }, 1 < p \le q < \infty; \eqno(2.15)
$$

\vspace{4mm}

$$
||F(u,v) ||_p  \le   ||u||_r \ ||v||_s, \ p,r,s > 1, \ 1/p = 1/r + 1/s. \eqno(2.16)
$$

 We denote for simplicity

 $$
 M'(d,r) = 0.5 M(d, d+r),
 $$
so that

$$
||\nabla w_t ||_r = M'(d,r) \cdot t^{ -1/2 - d(1-1/r)/2 }
$$

and
$$
 || \nabla  e^{t \Delta} v||_q \le  K_{BL}\left(d; r_0(p,q),p  \right) \cdot M'(d,r_0(p,q)) \
\cdot ||v||_p  \cdot t^{-(1+ d(1-1/r_0(p,q)) )/2 }, \eqno(2.17)
$$
where

$$
1 + \frac{1}{q} = \frac{1}{p} + \frac{1}{r_0(p,q)}. \eqno(2.18)
$$
 Note that $ M'(d,r) \le 2^{d-1}.  $ \par

\vspace{4mm}

\section{ Main result. }

\vspace{4mm}

{\it We suppose in this section that the initial function $ a = a(x) $ belong to the space $  L_d: \ ||a||_d < \infty,  $
and  such that } $  \Div a = 0. $ \par

\vspace{4mm}

{\bf 1.  Derivation of the basic inequalities.} \par

\vspace{4mm}

Let $  \delta $  be arbitrary fixed number from the set $ (0,1).  $  We denote

$$
 u_0 = u_0(x,t) = e^{t \Delta} a = [w_t*a](x) \eqno(3.1)
$$

and consider  together with T.Kato \cite{Kato1} the following  recursion:

$$
u_{n+1} = u_0 + G u_n, \ n=0,1,\ldots. \eqno(3.2)
$$
 Recall that the functions $ u_n = u_n(x,t)  $ are vector functions. \par
 We denote also  $ K_n = K_n(d,\delta) = K_n(a;d,\delta)  :=  $

 $$
  \sup_{ t > 0} || t^{(1-\delta)/2 } \ u_n  || BC([0,\infty); L_{d/\delta})  = \sup_{ t > 0} ||  t^{(1-\delta)/2  } \ u_n||_{d/\delta},
 \eqno(3.3)
 $$
\vspace{4mm}
$   K_n' = K'_n(d)= K'_n(a;d,\delta):=  $

$$
 \sup_{t > 0}  || t^{1/2} \  \nabla u_n  || BC([0,\infty); L_{d})  = \sup_{ t > 0} || t^{1/2} \ \nabla u_n||_{d}, \eqno(3.4)
 $$

\vspace{4mm}

$ K_n(T) = K_n(d,\delta;T) = K_n(a; d,\delta;T) := $

 $$
  || t^{(1-\delta)/2 } \ u_n  || BC([0,T); L_{d/\delta})  = \sup_{ t \in ( 0.T) } ||  t^{(1-\delta)/2  } \ u_n||_{d/\delta},
 \eqno(3.3')
 $$

\vspace{4mm}
$  K_n'(T) = K'_n(d;T) = K'_n(a; d;T):= $

$$
  || t^{1/2} \  \nabla u_n  || BC([0,T); L_{d})  = \sup_{ t \in ( 0,T)} || t^{1/2} \ \nabla u_n||_{d}, \eqno(3.4')
 $$
\vspace{4mm}
and correspondingly

$ K_0 = K_0(d,\delta) = K_0(a; d,\delta) := $

 $$
   || t^{(1-\delta)/2 } \ u_0  || BC([0,\infty); L_{d/\delta})  = \sup_{ t \ge 0} ||  t^{(1-\delta)/2  } \ u_0||_{d/\delta},
 \eqno(3.5)
 $$

\vspace{4mm}
$ K'_0 = K'_0(d,\delta) = K'_0(a; d,\delta) := $

$$
 K_0' = K'_0(d) =  || t^{1/2} \  \nabla u_0  || BC([0,\infty); L_{d})  = \sup_{ t \ge 0} || t^{1/2} \ \nabla u_0||_{d}. \eqno(3.6)
 $$

\vspace{4mm}
$ K_0(T) = K_0(d,\delta;T) = K_0(a; d,\delta;T) := $

$$
 K_0(T) = K_0(d,\delta;T) =  || t^{(1-\delta)/2 } \ u_0  || BC([0,T); L_{d/\delta})  = \sup_{ t \in(0,T) } ||  t^{(1-\delta)/2  } \ u_0||_{d/\delta},
 \eqno(3.5')
 $$

\vspace{4mm}

$ K'_0(T) = K'_0(d,\delta;T) = K'_0(a; d,\delta;T) := $
$$
 K_0'(T) = K'_0(d;T) =  || t^{1/2} \  \nabla u_0  || BC([0,T); L_{d})  = \sup_{ t \in (0,T) } || t^{1/2} \ \nabla u_0||_{d}. \eqno(3.6')
 $$

 Obviously,

 $$
  K'_n(d;T) < K'_n(d), \ \hspace{8mm}   K_n(d;T)  <  K_n(d).
 $$

 Moreover,

 $$
\lim_{T \to 0}  K_0(T) = 0, \hspace{8mm}  \lim_{T \to 0}  K_0'(T) = 0,
 $$
 see \cite{Kato1}. \par

\vspace{4mm}

 As we know,  see  (2.10),  (2.11),

$$
  K_0(d,\delta) \le  K_{BL}(d; d, d/(d-(1+ \delta))) \cdot M(d, d/(d-(1 + \delta )) )  \cdot ||a||_d, \eqno(3.7)
$$

$$
 K'_0(d)  \le
 0.5 \cdot K_{BL}(d; d,d) \cdot M \left(d, 1 \right)  \cdot ||a||_d. \eqno(3.8)
 $$

 Further, we find using (2.12) and (2.13):

 $$
 || G u_n||_{d/\delta} \le K_R(d/\delta) \cdot K_R(d) \cdot \int_0^t   || \nabla u_n(t-s)  ||_d \cdot ||u_n(s)||_{d/\delta} \ ds \le
 $$

$$
K_R(d/\delta) \cdot K_R(d) \cdot K_n \cdot K'_n \cdot \int_0^t (t-s)^{-1/2} \cdot s^{ - (1-\delta )/2  } \ ds  =
$$

$$
K_R(d/\delta) \cdot K_R(d) \cdot K_n \cdot K'_n \cdot t^{ -(1-\delta)/2 } \cdot \frac{\Gamma(1/2) \ \Gamma(\delta/2)}{\Gamma((1+\delta)/2)} =
$$

$$
K_R(d/\delta) \cdot K_R(d) \cdot K_n \cdot K'_n \cdot t^{ -(1-\delta)/2 } \cdot \frac{\sqrt{\pi} \ \Gamma(\delta/2)}{\Gamma((1+\delta)/2)}.
$$

 Therefore,

 $$
 K_{n+1} \le  K_{BL}(d; d, d/(d-(1+ \delta))) \cdot M(d, d/(d-(1 + \delta )) )  \cdot ||a||_d +
 $$

$$
K_R(d/\delta) \cdot K_R(d) \cdot K_n \cdot K'_n \cdot \frac{\sqrt{\pi} \ \Gamma(\delta/2)}{\Gamma((1+\delta)/2)}.
 \eqno(3.9)
$$
 We obtain analogously

$$
|| \nabla G u(t)||_d  \le K_R^2(d) \int_0^t (t-s)^{- ( 1+\delta )/2 } \ ||u(s)||_d \ ||\nabla u(s)||_d \ ds \le
$$

$$
K_R^2(d) \cdot K_n \cdot K'_n \cdot \int_0^t (t-s)^{-(1 + \delta)/2 } \ s^{ -1+\delta/2 } \ ds =
$$

$$
K_R^2(d) \cdot K_n \cdot K'_n \cdot t^{-( 1- \delta)/2} \int_0^1 (1-z)^{ -( 1+\delta)/2} \ z^{\delta/2-1}  \ dz =
$$

$$
K_R^2(d) \cdot K_n \cdot K'_n \cdot t^{-( 1- \delta)/2}  \cdot \frac{\Gamma ((1-\delta)/2) \ \Gamma(\delta/2)  }{\sqrt{\pi}}.
$$

 Following,

 $$
 K'_{n+1} \le  0.5 \cdot K_{BL}(d; d,d) \cdot M \left(d, 1 \right)  \cdot ||a||_d +
 $$

 $$
  K_R^2(d) \cdot K_n \cdot K'_n \cdot  \frac{\Gamma ((1-\delta)/2) \ \Gamma(\delta/2)  }{\sqrt{\pi}}.
 \eqno(3.10)
 $$

To sum up the local.  Let us denote

$$
S_1 = S_1(d,\delta) = K_{BL}(d; d, d/(d-(1+ \delta))) \cdot M(d, d/(d-(1 + \delta )) ),  \eqno(3.11)
$$

$$
J_1 = J_1(d,\delta) = K_R(d/\delta) \cdot K_R(d) \cdot  \frac{\sqrt{\pi} \ \Gamma(\delta/2)}{\Gamma((1+\delta)/2)}, \eqno(3.12)
$$

$$
S_2 =  S_2(d) =  0.5 \cdot K_{BL}(d; d,d) \cdot M \left(d, 1 \right), \eqno(3.13)
$$

$$
 J_2 = J_2(d,\delta) = K_R^2(d) \cdot  \frac{\Gamma ((1-\delta)/2) \ \Gamma(\delta/2)  }{\sqrt{\pi}}.
 \eqno(3.14)
$$
 We obtained the following system of recurrent inequalities for the vector sequence  $ (K_n, K'_n): $

$$
K_{n+1} \le K_0 + J_1 K_n K_n', \eqno(3.15)
$$
$$
K'_{n+1} \le K_0' + J_2 K_n K_n'  \eqno(3.16)
$$
 with initial conditions $ (K_0, K'_0),  $  where $  K_0 \le  S_1(d,\delta) \cdot ||a||_d, \  K'_0 \le  S_2(d) \cdot ||a||_d.   $  \par

\vspace{4mm}

{\bf 2. Auxiliary facts. } \par

\vspace{4mm}

 Let us consider the following non-linear recurrent inequality: $ x_n \ge 0,  $

 $$
 x_{n+1} \le \alpha + \beta x_n  + \gamma x^2_n, \ n = 0,1,2,\ldots \eqno(3.17)
 $$
with initial condition  $ x_0 = x(0) > 0. $ Denote

$$
D(\alpha, \beta, \gamma) = (\beta-1)^2 - 4 \alpha \gamma, \ Z(\alpha,\beta, \gamma) = \frac{1-\beta + \sqrt{D(\alpha,\beta,\gamma)}}{2 \gamma}.
\eqno(3.18)
$$

{\bf Lemma 1.} Let

$$
\alpha, \gamma > 0, \ \beta \ge 0, \ D(\alpha,\beta,\gamma) > 0, \  Z(\alpha,\beta,\gamma) > 0, \ x(0) < Z(\alpha,\beta,\gamma).
$$
Then

$$
\sup_n x_n \le Z(\alpha,\beta,\gamma). \eqno(3.19)
$$

 This assertion may be proved easily by means of induction over $ n. $ \par

\vspace{4mm}

{\bf Lemma 2.} Let us consider the following system of recurrent relations: $  x_n, y_n \ge 0,   $

$$
x_{n+1} \le \alpha_1 + \beta_1 x_n y_n, \eqno(3.20)
$$

$$
y_{n+1} \le \alpha_2 + \beta_2 x_n y_n, \ n = 0,1,2,\ldots. \eqno(3.21)
$$
with {\it positive}  initial conditions $ x_0 =  x(0), \ y_0 = y(0).  $ \par
 We retain last notations and add some news:

 $$
 Det_1 = Det_1(\alpha_1,\alpha_2, \beta_1, \beta_2) \stackrel{def}{=} \alpha_2 \beta_1 - \alpha_1 \beta_2,  \  Det_2 = - Det_1,\eqno(3.22)
 $$

 $$
  D_1 = D_1(\alpha_1,\alpha_2, \beta_1, \beta_2) \stackrel{def}{=}(Det_1 +1)^2 - 4 \alpha_1 \beta_2, \eqno(3.23)
 $$

 $$
 D_2 = D_2(\alpha_1,\alpha_2, \beta_1, \beta_2) \stackrel{def}{=}  (Det_2 +1)^2 - 4 \alpha_2 \beta_1. \eqno(3.24)
 $$
 Suppose $ \alpha_1, \alpha_2, \beta_1,\beta_2,  D_1, D_2 > 0, \   $

$$
Z(\alpha_1, Det_1, \beta_2) > 0, \hspace{5mm} Z(\alpha_2, Det_2, \beta_1) > 0, \
$$

$$
 0 < x(0) < Z(\alpha_1, Det_1, \beta_2),  \hspace{5mm}
0 < y(0) < Z(\alpha_2, Det_2, \beta_1).
$$
 We assert by induction based on the lemma 1:

 $$
 \sup_n x_n \le Z(\alpha_1, Det_1, \beta_2), \hspace{5mm} \sup_n y_n \le Z(\alpha_2, Det_2, \beta_1). \eqno(3.25)
 $$

\vspace{4mm}

{\bf 3. Main result. }

\vspace{4mm}
 We restrict ourselves in relations (3.15) - (3.16) only the time interval $  t \in (0,T), \ T < \infty. $ We obtain then the following
system of non-linear inequalities:

$$
K_{n+1}(T) \le K_0(T) + J_1 K_n(T) \cdot K_n'(T), \eqno(3.26)
$$
$$
K'_{n+1}(T) \le K_0'(T) + J_2 K_n(T) \cdot K_n'(T).  \eqno(3.27)
$$
 We have taken into account the notations (3.5') - (3.6'). It remains to use the estimates (3.25). Namely,
we denote using last notations

$$
s_1= s_1(d,\delta) = J_1(d,\delta)K_0'(T) -  J_2(d,\delta) K_0(T), \ s_2(d,\delta) =  - s_1(d,\delta); \eqno(3.28)
$$

$$
V_1(a(\cdot); d,\delta;T) = Z( K_0(T); s_1(d,\delta),  J_2(d,\delta)  ), \eqno(3.29)
$$

$$
V_2(a(\cdot); d,\delta; T) = Z( K_0'(T); s_2(d,\delta),  J_1(d,\delta)  ). \eqno(3.30)
$$
 We  have proved in fact the following assertion: \par

 \vspace{4mm}

 {\bf Theorem 3.1.} {\it Let as before } $ a \in L_d(R^d).  $
{\it  Define the set $ L = L(a; d,\delta) $ as a set  $ \{ T_0 \}  $ of all the
values $  T_0  $ as an arbitrary positive solutions of  inequalities}

 $$
 K_0(T_0) < V_1(a(\cdot); d,\delta, T_0), \eqno(3.31)
 $$

 $$
 K_0'(T_0) < V_2(a(\cdot), d,\delta; T_0). \eqno(3.32)
 $$
 {\it Then the lifespan of solution $ T* $ of NS equations is greatest than} $  T_0: $

 $$
    T* \ge \sup \{ T_0(a; d,\delta), \ t_0 \in L \} \stackrel{def}{=} \hat{T}  = \hat{T}(a; d,\delta) . \eqno(3.33)
 $$

 \vspace{4mm}

 {\bf Remark 3.1.}  As long as $ \lim_{T \to 0+} [K_0(T) + K_0'(T)] = 0 $ and

 $$
|s_1|, \ |s_2| \le  J_1(d,\delta)K_0'(T) +  J_2(d,\delta) K_0(T) \le J_1(d,\delta)K_0' +  J_2(d,\delta) K_0,
 $$
the system of inequalities (3.31), (3.32)  has at last one positive solution. \par

 \vspace{4mm}

 {\bf Remark 3.2.} Evidently,

 $$
    T \ge \sup_{\delta \in (0,1)} \hat{T}(a; d,\delta)  \stackrel{def}{=} \tilde{T}(a; d). \eqno(3.33)
 $$

\vspace{4mm}

\section{Simplification.}

\vspace{4mm}

 The system of inequalities (3.31) - (3.32) is very complicate. We aim to in this section simplification of this relations
in order to obtain more convenient   explicit view for the lower  bound for lifespan $ T. $  \par

\vspace{3mm}

{\bf 1.}  Note first of all that

$$
J_1 \le J^{(1)} = J^{(1)}(d,\delta) \stackrel{def}{=} \frac{9 d^2}{2 \delta^2}, \eqno(4.1)
$$

$$
 J_2 \le
J^{(2)} = J^{(2)}(d,\delta) \stackrel{def}{=}  \frac{81 d^2}{4 \sqrt{\pi} \delta (1-\delta)},  \eqno(4.2)
$$
as long as $  \delta \in (0,1). $\par

\vspace{3mm}

{\bf 2.} We denote

$$
 K^{(0)}(T) = K^{(0)}(a(\cdot); d,\delta; T) = \max( K_0(a(\cdot); d,\delta;T),  \ K'_0(a(\cdot); d,\delta;T) ),
$$

$$
X_n = X_n(T) = X_n(d,\delta;T) = \max( K_n(d,\delta;T), \ K_n(d,\delta;T) ),
$$

$$
J = J(d,\delta) = \max( J^{(1)}(d,\delta), \  J^{(2)}(d,\delta)),  \eqno(4.3)
$$

$$
\delta_0 := \frac{ 2 \sqrt{\pi}}{9 + 2 \sqrt{\pi}} \approx 0.282577,  \eqno(4.4)
$$

$$
\overline{J} = \overline{J}(d) = \min_{\delta \in (0,1)} \ J(d,\delta) = \frac{9 d^2}{2 \delta_0^2}= C_1 \cdot d^2,\eqno(4.5)
$$

$$
C_1 \approx 56.35566683; \eqno(4.5a)
$$
then

$$
 J(d,\delta) = J^{(1)}(d,\delta), \ \delta \in (0, \delta_0); \hspace{4mm} J(d,\delta) = J^{(2)}(d,\delta), \ \delta_0 \le \delta < 1.
$$

 It is a reason to name the pair of values $ (\delta_0, \overline{J})  $ as a critical value for considered problem. \par

\vspace{3mm}

{\bf 3.} For the variables $  X_n = X_n(T) $ we can write  one inequality:

$$
X_{n+1} \le K^{(0)}(T) + \overline{J} \ X^2_n  \eqno(4.6)
$$
with initial condition $ X_0 =  K^{(0)}(T). $  We conclude  applying lemma 1: \par

\vspace{4mm}

{\bf Theorem 4.1.}  {\it Define the value $  T_0  $ as follows: }

$$
\max( K_0(T_0), K_0'(T_0)) \le \frac{3}{16 \ \overline{J}} = C_2/d^2, \eqno(4.7)
$$

$$
C_2 \approx 0.0033270; \eqno(4.7a)
$$
{\it or equally}

$$
\max \left( \max_{t \in (0,T)} ||  t^{(1-\delta)/2  } \  w_t* a||_{d/\delta},
\max_{t \in (0,T)} ||  t^{1/2} \nabla w_t *a  ||_d  \right) \le \frac{3}{4 \ \overline{J}} = C_3/d^2. \eqno(4.7b)
$$

{\it  Then the lifespan $  T $  of solution NS system is greatest than $ T_0:   T \ge T_0. $  Moreover:}

$$
\max(K_n(T), K'_n(T)) \le \frac{3}{4 \ \overline{J}} = C_3/d^2, \eqno(4.8)
$$

$$
C_3 \approx 0.0133308333. \eqno(4.8a)
$$

\vspace{4mm}
  We have for instance in the ordinary three-dimensional case $  d=3 $

$$
C_2/d^2 \approx 0.00036967, \hspace{6mm} C_3/d^2  \approx 0.0014767.
$$

\vspace{4mm}

{\bf 4. } In order to use the theorem 4.1 we need to derive a simple  estimate  for the values $ K_0(T), \ K_0'(T) $
as $ T \to 0+. $ \par

\vspace{3mm}

{\bf $ \alpha. $ Estimation of } $  K_0(T). $ \par
 Suppose in addition to the condition $ ||a||_d < \infty  $ that for some \\
 $ \theta \in (0, \min(1, (d-1)/\delta) $

$$
||a||_{d + \theta} < \infty. \eqno(4.9)
$$

 Repeating the consideration for the inequality (2.10)  and taking into account the inequalities $ K_{BL}(\cdot) \le 1, \
 M(d,r) \le 2^d $  we deduce:

 $$
 ||u_0||_{d/\delta} \le t^{-( 1-\delta)/2   } \cdot t^{(\theta \delta)/(2d)  } \cdot 2^{d + \theta} \cdot ||a||_{d + \theta};
 $$
 therefore

 $$
 K_0(T) \le T^{ \frac{\theta \delta}{2d} } \cdot 2^{d + \theta}  \cdot ||a||_{d + \theta}. \eqno(4.10)
 $$

\vspace{3mm}

{\bf $ \beta. $ Estimation of } $  K'_0(T). $  \par

 We demonstrate in this pilcrow  a different method. Namely, assume in addition that the initial condition
$ a(\cdot)  $  belongs to the Sobolev space $ W_1^d(R^d), $  which consists on all the (measurable) functions $ a: R^d \to R^d $
which finite semi-norm

$$
||a||W_1^d(R^d) = || \nabla a||_d < \infty. \eqno(4.11)
$$
 We get using again Young's inequality:

 $$
 K'_0(T) = \sup_{ t \in (0,T)} \left[ t^{1/2} \ \nabla (w_t)*a \right]_d = \sup_{ t \in (0,T)} \left[t^{1/2} \ w_t* \nabla a \right]_d \le
 $$

$$
\sup_{t \in (0,T)} \left[ t^{1/2} \ ||w_t||_1 \ ||\nabla a||_d    \right] = \sup_{t \in (0,T)} \left[ t^{1/2} \ ||\nabla a||_d    \right] =
$$

$$
\sqrt{T} \cdot ||a||W_1^d(R^d), \eqno(4.12)
$$
since $ ||w_t||_1  = 1.  $\par

\vspace{4mm}

{\bf Remark 4.1.} The variables $  K_0(T), \ K_0'(T) $ dependent in particular on the initial condition $ a(\cdot): $

$$
K_0(T)=  \ K_0(T; a(\cdot))); \hspace{5mm} K_0'(T)= K_0'(T; a(\cdot)).
$$
If  it is so little that

$$
\max( K_0(a(\cdot)), K_0'(a(\cdot))) \le \frac{3}{16 \ \overline{J}} = C_2/d^2, \eqno(4.13)
$$
then  we can choose $ T = \infty,  $ i.e. this solution $ u = u(x,t) $ is {\it  global.}\par

 Recall that

$$
 K_0(a(\cdot)) = K_0(a(\cdot); \infty); \hspace{6mm} K_0'(a(\cdot)) = K_0'(a(\cdot); \infty).
$$

\vspace{4mm}

\section{Mixed norm estimates for  solution.}

\vspace{4mm}
 {\it We suppose during this section } $  ||a||_d <\infty.  $ \par

\vspace{4mm}

 It is known, see \cite{Cui1}, \cite{Kato1}-  \cite{Kato2} that  the global in time solution $ u(x,t) = u(t)  $ obeys the
property

$$
\lim_{t \to \infty} ||u(t)||_q = 0,  \ q > d.
$$
( "Energy" decay). \par
 The case  $  q = 2  $  was investigated in  \cite{Ogawa1}; see also reference therein. \par
 We  want  clarify in this section this fact; i.e. give the quantitative estimates one of main result if the article
 \cite{Kato1}. \par

Recall that the so-called {\it mixed,  } or equally {\it anisotropic } $ (p_1, p_2) $   norm  $ ||u||^*_{p_1,p_2} $ for the function of "two"
variables $ u = u(x,t), \ x \in R^d, \ t \in R^1_+ $  is defined as follows:

$$
||u||^*_{p_1,p_2} = \left( \int_{R^d} \left[ \int_0^{\infty} |u(x,t)|^{p_1} \ dx  \right]^{p_2/p_1} \ dt \right)^{1/p_2}
$$
with evident modification in the case when $ p_2 = \infty: $

$$
||u||^*_{p_1,\infty} = \sup_{t \in (0,T)} \int_{R^d} \left[ \int_0^{\infty} |u(x,t)|^{p_1} \ dx  \right]^{1/p_1}.
$$

 We introduce here the following {\it weight} mixed  norm, more precisely, the family of norms as follows:

$$
|||u|||^*_{q;T} = \sup_{t \in (0,T)} \left[ t^{(1-d/q)/2 } \ ||u(\cdot,t)||_q \right], \ q \ge d. \eqno(5.1)
$$

 Let us introduce  some new notations.

 $$
 \theta_1 = \frac{dq}{d(q+1) -q(\delta+1) }, \ \theta_2 = \frac{q}{\delta+1},     \ \theta_3 = \frac{d}{\delta}, \ \theta_4 = d;
 $$

 $$
 \psi(d, q,\delta) =  \psi(a;d.q,\delta) = K_{BL}(d; \theta_1, \theta_2) \  K^0(a; d,q,\delta; T) \ K^0(a; d,q,\delta; T)' \times
 $$

 $$
 K_R(d/\delta) \ K_R(d) \ M(d, \theta_1) \ B \left(\frac{1-\delta}{2} + \frac{d}{2q}, \ \frac{\delta}{2}  \right) +
 $$

 $$
 0.5 \cdot K_{BL}(d; d,d) \cdot M \left(d, \frac{d^2}{d - 1} \right)  \cdot ||a||_d,\eqno(5.2)
 $$
 where $ B(\cdot, \cdot) $ is the classical Beta-function;

 $$
 \psi(q) = \psi(a;q) =  \inf_{\delta \in (0,1)}  \psi(a;d.q,\delta). \eqno(5.3)
 $$

 Note that $ 1 < \theta_j < \infty, \ j = 1,2,3,4.  $ \par

\vspace{4mm}

{\bf Theorem 5.1.} {\it  Let the lifespan of solution of NS system }  $ T  $ {\it  be positive; may be infinite. Then  }

$$
|||u|||^*_{q;T} \le \psi(q), \ q \ge d. \eqno(5.4)
$$

\vspace{4mm}

{\bf Proof.} We follow T.Kato \cite{Kato1}. Indeed,

$$
- G u(t) = \int_0^t w_{t-s}(\cdot)*F(u(\cdot),s) \ ds.
$$
 We use the triangle inequality for the $ L_q(R^d) $ norm:

 $$
 ||Gu(t)||_q \le \int_0^t ||  w_{t-s}(\cdot)*F(u(\cdot),s)||_q \ ds. \eqno(5.5)
 $$
 The Young inequality give us

 $$
  ||  w_{t-s}(\cdot)*F(u(\cdot),s)||_q \le K_{BL}(d; \theta_1, \theta_2) \  K_{R}(\theta_3) \ K_{R}(\theta_4) ||  w_{t-s}(\cdot) ||_{\theta_1} \  ||F||_{\theta_2} \le
 $$

  $$
  K_{BL}(d; \theta_1, \theta_2)  \cdot M(d, \theta_1) \cdot (t-s)^{-d(1 - 1/\theta_1)/2} \cdot K_{R}(\theta_3) \cdot K_{R}(\theta_4) \times
 $$

 $$
    ||u(\cdot,s)||_{\theta_3}  \cdot  ||\nabla u(\cdot,s)||_{\theta_4}, \eqno(5.6)
 $$
 we have taken into account the norm of Riesz transform. Further, if $ t \in (0,T) $ then

$$
 ||u(\cdot,s)||_{d/\delta} \le   K^0(a; d,q,\delta; T) \ s^{ -(1-\delta)/2 },
$$

$$
 || \nabla u(\cdot,s)||_d \le   K^0(a; d,q,\delta; T)' \ s^{ -1/2 }.
$$

 We conclude substituting into the inequality (5.5):

$$
||Gu(t)||_q \le C(a;d,q,\delta,T) \cdot  \int_0^t (t-s)^{ - ( 1 + \delta - d/q)/2  } \ s^{-1 + \delta/2 } \ ds =
$$

$$
C(a;d,q,\delta,T) \ t^{ -(1 - d/q)/2  } \ B \left(\frac{1-\delta}{2} + \frac{d}{2q}, \ \frac{\delta}{2}  \right). \eqno(5.7)
$$

 As long as

 $$
 ||u(t)||_q \le ||a||_q + ||Gu(t)||_q
 $$
and the value $ ||a||_q $ was estimated in (2.11), we deduce after simple calculations

$$
||u(t)||_q \le \ t^{ -(1 - d/q)/2  } \ \psi(a; d,q, \delta,T)
$$
or equally

$$
|||u|||^*_{q;T} \le \psi(a; d,q, \delta,T). \eqno(5.8)
$$
 It remains to take the minimum over $ \delta; \ \delta \in (0,1).  $ \par

\vspace{4mm}

{\bf Remark 5.1.} If we define the so-called mixed Grand Lebesgue norm  $ |||u|||^*G(\psi; T)   $ as follows:

$$
|||u|||^*G(\psi; T):=  \sup_{q \ge d} \left[ \frac{|||u|||^*_{q;T}}{\psi(a; q)}\right],  \eqno(5.9)
$$
then the assertion of the theorem 5.1 may be rewritten as follows: under the conditions of theorem  5.1

$$
|||u|||^*G(\psi(a; T)) \le 1.  \eqno(5.10)
$$
 The  detail investigation with applications of these norm and correspondent spaces  see, e.g. in \cite{Fiorenza3} - \cite{Jawerth1},
\cite{Kozachenko1}, \cite{Liflyand1}, \cite{Ostrovsky2}, \cite{Ostrovsky107}. \par

 \vspace{4mm}

 Let us investigate here the following {\it weight} mixed  norm with derivative, (Sobolev's weight norm) for solution $ u = u(x,t),  $
more precisely, the family of {\it semi-norms} as follows:

$$
|||u|||^{**}_{q;T} = \sup_{t \in (0,T)} \left[ t^{(1-d/2q) } \ || \nabla u(\cdot,t)||_q \right], \ q \ge d. \eqno(5.11)
$$

\vspace{4mm}

 New notations and restrictions:

 $$
 1 + \frac{1}{q} = \frac{1}{\theta_5} + \frac{1}{d}, \ 1 + \frac{1}{q} = \frac{1}{\theta_6} + \frac{1}{\theta_7},
 $$

$$
\frac{1}{\theta_7} = \frac{\delta}{d} + \frac{1}{d},
$$

$$
\nu(a;d,\delta, q, T) =  K_{BL}(d; q,\theta_5) \cdot M'(d, \theta_5) ||a||_d +
$$

$$
K_{BL}(d; \theta_6, \theta_7) \ K_R(\theta_6) \  K_R(\theta_7) \ K(a; d,\delta, T) \ K'(a; d,\delta, T) \ M'(d,\theta_6),
$$

$$
\nu(q) = \nu(a; q):=  \inf_{ \delta \in (0,1)}  \nu(a;d,\delta, q, T).
$$

\vspace{4mm}

{\bf Theorem 5.2.} {\it  Let the lifespan of solution of NS system }  $ T  $ {\it  be positive; may be infinite. Then  }

$$
|||u|||^{**}_{q;T} \le \nu(q), \ q \ge d. \eqno(5.12)
$$

\vspace{4mm}
 {\bf Proof.} First of all we estimate the influence  of the initial condition $ a = a(x).  $  Namely,

 $$
 ||\nabla u_0||_q  =  || \nabla w_t * a ||_q \le K_{BL}(d; q,\theta_5) \cdot || \nabla w_t ||_{\theta_5} \cdot
 ||a||_d =
 $$

$$
t^{-(1 - d/(2q))} \cdot  K_{BL}(d; q,\theta_5) \cdot M'(d, \theta_5) ||a||_d. \eqno(5.13)
$$
  Further,

  $$
  ||\nabla Gu(t)||_q = ||\int_0^t \nabla w_{t-s} * F(u(s)) \ ds ||_q \le \int_0^t ||\nabla w_{t-s} * F(u(s)) ||_q \ ds \le
  $$

$$
K_{BL}(d; \theta_6, \theta_7) \ K_R(\theta_6) \  K_R(\theta_7)  \int_0^t ||\nabla w_{t-s}||_{\theta_6} \ ||F(u(s))||_{\theta_7} \ ds \le
$$

$$
K_{BL}(d; \theta_6, \theta_7) \ K_R(\theta_6) \  K_R(\theta_7)  \int_0^t ||\nabla w_{t-s}||_{\theta_6} \ ||u(s)||_{d/\delta} \
||\nabla u(s) ||_{d} \ ds. \eqno(5.14)
$$
 As long as

 $$
 ||u(s)||_{d/\delta} \le s^{-( 1-\delta )/2 } K(a; d,\delta, T), \ ||\nabla u||_d \le K'(a; d,\delta, T) s^{-1/2}, \eqno(5.15)
 $$

 $$
 ||\nabla w_{t-s}||_{\theta_6} = M'(d,\theta_6) \ t^{-1/2 - d(1-1/\theta_6)/2 }, \eqno(5.16)
 $$
we obtain  substituting into (5.14):  $ ||\nabla Gu(t)||_q  \le $

$$
K_{BL}(d; \theta_6, \theta_7) \ K_R(\theta_6) \  K_R(\theta_7) \ K(a; d,\delta, T) \ K'(a; d,\delta, T) \ M'(d,\theta_6) \times
$$

$$
\int_0^t  (t - s)^{-1/2 - d(1-1/\theta_6)/2 } \ s^{ (-1 + \delta/2)} \ ds  = t^{-( 1 - d/(2q))} \cdot
B \left( 0.5 - d(1-\theta_6)/2, \ \delta/2   \right) \times
$$

$$
K_{BL}(d; \theta_6, \theta_7) \ K_R(\theta_6) \  K_R(\theta_7) \ K(a; d,\delta, T) \ K'(a; d,\delta, T) \ M'(d,\theta_6). \eqno(5.17)
$$
We  get summing (5.13) and (5.17):

$$
||u(t)||_q \le t^{-( 1 - d/(2q))} \cdot \nu(a;d,\delta, T)
$$
and after minimization over $ \delta $

$$
\sup_{t \in (0,T)} \left[t^{( 1 - d/(2q))} \ ||u(t)||_q \right]  \le \nu(q), \eqno(5.18)
$$
Q.E.D. \par

\vspace{4mm}

 \section{Concluding remarks. }

 \vspace{4mm}

{\bf 1.}  It is known \cite{Giga1}, \cite{Giga2},  \cite{Kato1}, \cite{Kato2}  etc. that in general case, i.e.
when the value $ \epsilon = ||a||_d $ is not sufficiently small, then the lifespan of solution of NS equation $ T $
may be finite (short-time solution).  Perhaps, it is  self-contained  interest to find a quantitative computation of the {\it exact
value} $ T. $ \par
  For the non-linear Schr\"odinger's equation the estimate

 $$
 T \ge  \exp (C /\epsilon)
 $$
was obtained in the recent article \cite{Ikeda1}.\par

\vspace{4mm}
{\bf 2.} At the same considerations may be provided for the NS equations with external force $  f = f(x,t): $

$$
\partial{u}_t = \Delta u + (u \cdot \nabla)u  +   Q \cdot \nabla \cdot (u \otimes u) + f(x,t), \ x \in R^d, \ t > 0; \eqno(6.1)
$$

$$
u(x,0) = a(x), \ x \in R^d.
$$
see  \cite{Giga1} -  \cite{Giga4}, \cite{Koch1},  \cite{Masuda1}, \cite{Ogawa1}, \cite{Temam1}.\par
 More detail, the considered here problem may be rewritten as follows:

$$
u(x,t) = e^{t \Delta} a(x)  + G [u](t) \stackrel{def}{=} u_0(x,t) + G [u](t) + v[f](x,t), \eqno(6.2)
$$
where

$$
v[f](x,t) = v(x,t) = v = \int_0^t ds \int_{R^d} w_{t-s}(x-y) \ f(y,s) \ dy =
$$

$$
\int_0^t  w_{t-s}(\cdot)*f(\cdot,s) \ ds. \eqno(6.3)
$$

 In order to formulate a new result we introduce new Banach spaces on the (measurable) functions of two variables $ f(x,t), \ x \in R^d, \ t > 0:  $

$$
|||f|||_{\theta, \lambda} := \sup_{ s > 0} \left[ \frac{||f(\cdot,s)||_{\theta}}{ s^{\lambda}} \right], \eqno(6.4)
$$

$ \theta = \const \ge 1, \ \lambda = \const \in (-1,0).  $  The space of all the functions $ \{f(\cdot,\cdot) \}  $ with
finite such a norm will be denoted $ A(\theta,\lambda):  $

$$
 A(\theta,\lambda) = \{f: |||f|||_{\theta, \lambda} < \infty  \}. \eqno(6.5)
$$
 {\it We assume in this subsection}

$$
f \in A(\theta_1, \lambda_1) \cap  A(\theta_2, \lambda_2). \eqno(6.6)
$$

 Define  for arbitrary $  \delta \in (0,1) $ the values $ r_1, r_2 $ as follows:

$$
1 + \frac{\delta}{d} = \frac{1}{r_1} + \frac{1}{\theta_1}, \hspace{4mm} 1 + \frac{1}{d} = \frac{1}{r_2} + \frac{1}{\theta_2}
$$
with the following restrictions:

$$
r_{1,2} > 1, \theta_{1,2} \ge 1,  \hspace{4mm} -1 < \lambda_{1,2} < 0,
$$

$$
d \left(1-\frac{1}{r_1} \right) < 2,  \hspace{4mm}    d \left(1-\frac{1}{r_2} \right) < 1, \eqno(6.7)
$$

$$
\frac{d}{2} \left( 1 - \frac{1}{r_1} \right) - 1 - \lambda_1  = \frac{1-\delta}{2}, \hspace{4mm}
\frac{d}{2} \left( 1 - \frac{1}{r_2} \right)-1 - \lambda_2  = \frac{1}{2}.
$$

 We estimate using once more Young's inequality:

 $$
 ||v ||_{d/\delta} \le \int_0^t || w_{t-s}*f(s)||_{d/\delta} \ ds  \le K_{BL}(r_1, \theta_1) \int_0^t ||w_{t-s}||_{r_1} \cdot
 ||f(s)||_{\theta_1} ds \le
 $$

$$
K_{BL}(r_1, \theta_1) \ M(d, r_1) \ ||f||_{\theta_1, \lambda_1} \ \int_0^t (t-s)^{-d( 1-1/r_1 ) } \ s^{\lambda} \ ds =
$$

$$
K_{BL}(r_1, \theta_1) \ M(d, r_1) \ ||f||_{\theta_1, \lambda_1} \ t^{1 - 0.5 d (1-1/r_1) + \lambda_1 }
B(1-d( 1-1/r_1 ), 1 + \lambda_1) =
$$

$$
 t^{ -(1-\delta)/2 } \  K_{BL}(r_1, \theta_1) \ M(d, r_1) \ ||f||_{\theta_1, \lambda_1} \
B(1-d( 1-1/r_1 ), 1 + \lambda_1). \eqno(6.8)
$$

 \vspace{4mm}

  We find analogously

$$
\nabla v = \int_0^t ds \int_{R^d} \nabla w_{t-s}(x-y) \ f(y,s) \ ds;
$$

$$
||\nabla v||_d \le K_{BL}(r_2, \theta_2) \ M'(d,r_2) \ ||f||_{\theta_2, \lambda_2} \ \int_9^t (t-s)^{-1/2 - 0.5 d (1-1/r_2)  } \
s^{\lambda_2} \ ds =
$$

$$
 K_{BL}(r_2, \theta_2) \ M'(d,r_2) \ ||f||_{\theta_2, \lambda_2} \ B( 1/2 - 0.5 d (1-1/r_2), \lambda_2 + 1  ) \ t^{-1/2}. \eqno(6.9)
$$

Summing with the estimations for $ ||u_0||_{d/\delta}, \ ||\nabla u_0||_{d}  $  with corespondent estimations  (6.8), (6.9) for
$ ||v||_{d/\delta}, \ ||\nabla v||_d $ we conclude that if the norms
$ ||a||_d, \ ||f||_{\theta_1, \lambda_1}, \    ||f||_{\theta_2, \lambda_2} $ are finite. then $ T > 0;  $
if these norms are sufficiently small,  then $  T = \infty.  $ \par

\vspace{4mm}

 {\bf 3.} Analogously to the content of this report may be considered a more general case of abstract (linear or not linear) parabolic
 equation of a view

 $$
 \partial{u}_t  = A u + F(u, \nabla u; x,t)  + f(x,t), \hspace{5mm} u(x,0) = a(x).
 $$
 The detail investigation of this case when the initial condition and external force belong to some
 Sobolev's space may be found, e.g. in \cite{Taylor1} -  \cite{Taylor3}, \cite{Iwashita1}, \cite{Solonnikov1}. \par

\vspace{4mm}

{\bf 4.}  Let us consider the general non-linear parabolic equation  (may be multivariate, i.e. system of equations) of a view

$$
\frac{\partial u}{\partial t} = Au + F(u, \nabla u)
$$
with initial condition $ u(x,0+) = a(x). $ Here $  A  $ be negative definite linear operator, may be unbounded, for example

$$
A u = \sum_{i,j = 1}^d b_{i,j}(x) \frac{\partial^2 u}{\partial x_i \ \partial x_j},
$$
i.e $  A  $ be strictly elliptical differential operator of a second order with bounded coefficients:

$$
0 < \inf_{\xi: ||\xi|| = 1} \inf_{x \in R^d} \sum_{i,j=1}^d b_{i,j}(x) \xi_i \xi_j  \le
\sup_{\xi: ||\xi|| = 1} \sup_{x \in R^d} \sum_{i,j=1}^d b_{i,j}(x) \xi_i \xi_j  < \infty,
$$
satisfying the H\"older's condition: for some positive $ \alpha = \const \in (0,1]  $

$$
\max_{i,j}|b_{i,j}(x)  - a_{i,j}(y)| \le  C \cdot ||x - y||^{\alpha}.
$$

 We set ourselves a goal to obtain a positive lower  estimate for lifespan $ T $ for solution likewise to the
third section, following, e.g. \cite{Taylor3}, p. 272-275. \par
We understand as before in the capacity of solution  the mild solution:

$$
u(t) = e^{A t} a + \int_0^t e^{(t-s)A} \ F(u(s), \nabla u(s)) \ ds \stackrel{def}{=}
$$

$$
 e^{A t} a + \int_0^t e^{(t-s)A}  \ \Phi(u(s)) \ ds \stackrel{def}{=} \Psi[u](t),
$$
where $ \Phi(u(s) = F(u(s), \nabla u(s)). $ \par
\vspace{3mm}
We enumerate the conditions. \par
\vspace{3mm}
{\bf A.}  There exist two Banach spaces $ X, \ Y $ such that $ a(\cdot) \in X  $ and such that the semigroup
$  \{ S_t\} = \{ \exp(At)  \}, \ X \to X \ t \ge 0  $ is strong continuous.  \par

\vspace{3mm}
{\bf B.} The operator  (non linear, in general case) $ \Phi: X \to Y  $ is locally Lipshitz map. \par

\vspace{3mm}
{\bf C.}

$$
\forall t > 0 \hspace{4mm}   S_t = e^{A t}: Y \to X
$$
and  moreover

$$
\exists \gamma \in (0,1), \  \exists C = C(\gamma) \ \Rightarrow \ || e^{A t} ||(Y \to X) \le C(\gamma) \ t^{ - \gamma }, \
t \in (0,1).
$$
 Many examples of such a semigroups are given in \cite{Taylor3}, p. 274 \ - \ 286. \par

\vspace{4mm}
 Let $ \alpha = \const \in (0,1); $  we define

 $$
 Z = Z(\alpha) = \{ u \in C([0,T_1],X: \ u(0) = a, \ ||u(t) - f||X \le \alpha   \}. \eqno(6.10)
 $$
 There exists a value $  T_2  = T_2(\alpha)> 0 $ for which
$$
\sup_{t \in (0,T_2)} || e^{tA} a - a||X \le \frac{\alpha}{2}.
\eqno(6.11)
$$
 Denote also

 $$
 K_1 = K_1(\alpha) = \sup_{u \in Z(\alpha)} || \Phi(u)(\cdot) ||Y < \infty;
 $$
then

$$
|| \int_0^t e^{ A(t-s)} \Phi(u(s)) \ ds ||X \le K_1(\alpha) \ C(\gamma) \frac{t^{1-\gamma}}{1-\gamma}. \eqno(6.12)
$$

We can choose  $ T_3 $ such that

$$
K_1(\alpha) \ C(\gamma) \frac{T_3^{1-\gamma}}{1-\gamma} < \frac{\alpha}{2}.
$$

 As long as

 $$
 \forall u,v \in Z(\alpha) \ \Rightarrow || \Phi(u)(s) - \Phi(v(s))||Y \le K_2(\alpha) ||u(s) - v(s)||X, \
 K_2 = K_2(\alpha) < \infty,
 $$
we deduce

$$
||\Psi[u](t) - \Psi[v](t)|| \le \int_0^t e^{A s } || \Phi(u)(s) - \Phi(v(s))||Y  \ ds \le
$$

$$
K_2(\alpha) C(\gamma) \frac{ t^{1-\gamma} }{1-\gamma} \sup_{s \in (0,T_4)} ||u(s) - v(s) ||X \le
$$

$$
0.5  \sup_{s \in (0,T_4)} ||u(s) - v(s) ||X,
$$
if

$$
K_2(\alpha) C(\gamma) \frac{ T_4^{1-\gamma} }{1-\gamma} \le 0.5.
$$
Thus, we can take  $ T = \min(T_1, T_2, T_3, T_4). $ \par

 \vspace{4mm}

On the other hand, the examples of blow up in finite time solutions  of non-linear parabolic equations see, e.g. in articles
\cite{Ball1}, \cite{Marino1}. Examples of non-uniqueness for this equations see in \cite{Ladyzhenskaya1}.  \par

\vspace{4mm}

{\bf 5.} It may be investigated in addition analogously the {\it  boundary } value problem, for instance, Dirichlet
or Neuman, when the variable $ x $  belongs to some  domain $  \Omega $ (bounded or not)   with  boundary of the class
$  C^{1,1}. $\par

\end{document}